\documentclass[oneside,a4paper]{amsart}
\usepackage{quad}

\begin{document}

\title{A quest for 5-point condition a la Alexandrov}
\author{Anton Petrunin}
\date{}
\begin{abstract}
Here I give a description of Alexandrov 4-point comparison
via quadratic forms
and then propose a natural 5-point condition which might have some future.
Consider this note as a letter from me --- do not take it seriously.
\end{abstract}
\maketitle

\section*{Associated form}

Let us construct a quadratic form $W_{\bm{x}}$ on $\RR^{n-1}$ 
for given $n$-point array $\bm{x}\z=(x_1,\dots,x_n)$ in a metric space $X$.

Fix a nondegenerate simplex $\triangle$ in $\RR^{n-1}$;
denote by $v_1,\dots,v_n$ its vertices;
if $e_i$ is the standard basis on $\RR^{n-1}$,
we may assume that $v_i=e_i$ for $i<n$ and $v_n=0$.

Let us denote by $|a-b|_X$ the distance between points $a$ and $b$ in the metric space $X$.
Set
\[W_{\bm{x}}(v_i-v_j)=|x_i-x_j|^2_X\] 
for all $i$ and $j$.
Note that this identity define $W_{\bm{x}}$ uniquely.

The constructed quadratic form $W_{\bm{x}}$ will be called \emph{associated form} for the point array $\bm{x}$ with respect to the simplex $\triangle$.

Note that an array $\bm{x}=(x_1,\dots,x_n)$ in a metric space $X$ is isometric to an array in Euclidean space if and only if 
$W_{\bm{x}}(v)\ge 0$
for any $v\in \RR^{n-1}$.

In particular, for $n=3$ the
condition $W_{\bm{x}}\ge 0$ means that for the points $x_1$, $x_2$ and $x_3$,
all three triangle inequalities hold. 

\section*{Alexandrov's 4-point comparison}

Now let us discuss the relation between form $W_{\bm{x}}$
and geometry of the space for $n=4$.
In this case $\triangle$ is a tetrahedron on $\RR^3$.

\begin{wrapfigure}{r}{52mm}
\begin{lpic}[t(-5mm),b(3mm),r(0mm),l(0mm)]{pics/quad(1)}
\lbl[]{11,-2;$\quadra(4)$}
\lbl[]{40,-2;$\quadra(3)$}
\end{lpic}
\end{wrapfigure}

From the case $n=3$ 
it follows that $W_{\bm{x}}$ 
is nonnegative on every plane parallel to a face of the tetrahedron $\triangle$.
In particular, $W_{\bm{x}}$ can have at most one negative eigenvalue.

Assume $W_{\bm{x}}(w)<0$ for some $w\in\RR^3$.
It follows from above that
$w$ is transversal to each of 4 planes parallel to a faces of $\triangle$.
Let us project $\triangle$ along $w$ to a transversal plane. 
Note that after the projection the 4 vertices of $\triangle$ lie in general position; 
i.e., no three of them lie on one line.
Therefore in the projection we can see one of two combinatorial pictures shown on the diagram.
It is easy to see that the combinatorics of the picture does not depend on the choice of $w$.

If we see the diagram on the left we say that $\bm{x}$ is 
of type $\quadra(4)$ and otherwise we say that it is of type $\quadra(3)$.

The following statements give a connection between the associated forms $W_{\bm{x}}$
and the curvature bounds in the sense of Alexandrov.
Their proofs are left to the reader.

Assume $X$ is a complete space with intrinsic metric.
Then
\begin{itemize}
\item If $W_{\bm{x}}\ge 0$ 
for any quadrilateral $\bm{x}=(x_1,\dots,x_4)$ 
then $X$ is isometric to a closed convex set in a Hilbert space. 
\item $X$ has no quadruples of type $\quadra(3)$ if and only if 
$X$ has nonnegative curvature in the sense of Alexandrov.
\item $X$ has no quadruples of type $\quadra(4)$ if and only if 
$X$ is a $\CAT[0]$ space.
\end{itemize}

\section*{5-point conditions}

Let us try to do the same 
for 5-points arrays $\bm{x}=(x^1,\dots,x^5)$ in a metric space.
In this case the associated form $W_{\bm{x}}$ is defined on $\RR^4$
and it has to be nonnegative on any plane which is parallel to any of 10 two-dimensional faces of the $4$-simplex $\triangle\subset \RR^4$.
In particular $W_{\bm{x}}$ has at most two negative eigenvalues.

In the case if $W_{\bm{x}}$ has exactly two negative eigenvalues,
one can choose a plane $\Pi$ 
such that the restriction of $W_{\bm{x}}$ to $\Pi$ is negative.
Let us project $\triangle$ along $\Pi$ to a transversal plane.
The same argument as in case $n=4$ shows that after projection
the vertices of $\triangle$  lie in general position.
Therefore we may get one of the following three combinatorial pictures.

\begin{wrapfigure}{r}{70mm}
\begin{lpic}[t(-8mm),b(3mm),r(0mm),l(0mm)]{pics/penta(1)}
\lbl[]{10,-3;$\penta(5)$}
\lbl[]{34,-3;$\penta(4)$}
\lbl[]{57.5,-3;$\penta(3)$}
\end{lpic}
\end{wrapfigure}
I.e., for any 5-point array $\bm{x}$
either $W_{\bm{x}}$ has at most one negative eigenvalue
or it has exactly two negative eigenvalues and belongs to one of these three types $\penta(5)$, $\penta(4)$ or $\penta(3)$.

We may consider the metric spaces
which do not have a 5-point arrays of some of these types.
For example, 
denote by $\penta(\hat 3,\hat 4)$ the class of all complete length-metric spaces 
without 5-point arrays of type $\penta(3)$ and $\penta(4)$.

Let us list some easy observations about these new classes of metric spaces.

\begin{enumerate}[(i)]
\item\label{i} Any $\CAT[0]$ space is a $\penta(\hat 3,\hat 4,\hat 5)$ space.
In other words,
the form $W_{\bm{x}}$ for any 5-point array $\bm{x}$ in any $\CAT[0]$ space 
has at most one negative eigenvalue.
\item\label{ii} Any space with nonnegative curvature in the sense of Alexandrov has no  5-point arrays of type $\penta(3)$ and $\penta(4)$; 
i.e., it is a $\penta(\hat 3,\hat 4)$ space.
\item\label{iii} If a complete Riemannian manifold $M$ has no
5-point array of type $\penta(5)$ then it is a simply connected and it has   non-positive sectional curvature.
In other words, if a complete Riemannian manifold $M$ belongs to $\penta(\hat 5)$ then it is $\CAT[0]$.
\end{enumerate}

From the observations (\ref{i}) and (\ref{ii}) it follows that the $\penta(\hat 3,\hat 4)$
includes all $\CAT[0]$ spaces 
as well as the spaces with nonnegative curvature in the sense of Alexandrov.

\begin{thm}{Question}
 Do $\penta(\hat 3,\hat 4)$ spaces have meaningful geometry?
\end{thm}
In particular we do not know the answer to the following question.

\begin{thm}{Question}
Let $X$ be a complete $\penta(\hat 3,\hat 4)$ space with intrinsic metric.
Assume two minimizing geodesics in $X$ share three distinct points.
Does it imply that these geodesics share an arc containing these three points?
\end{thm}

\begin{thm}{Question}
Let $X$ be the coordinate plane with metric induced by norm.
Assume $X$ is a $\penta(\hat 3,\hat 4)$ space.
Is it true that $X$ has to be isometric to Euclidean plane?
\end{thm}

Here is an other question related to the observations (\ref{i}) and (\ref{iii}).

\begin{thm}{Question}
Is it true that any complete  $\penta(\hat 3,\hat 4,\hat 5)$ space with intrinsic metric has to be $\CAT[0]$?
\end{thm} 

\end{document}